\documentclass[preprint,3p,sort]{elsarticle}

\usepackage{listings}
\usepackage{amsfonts}
\usepackage{amsthm}
\usepackage{graphicx}
\usepackage{epstopdf}
\usepackage{algorithmic}
\usepackage{amsmath}
\usepackage{amssymb}
\usepackage{hyperref}
\usepackage[T1]{fontenc}
\usepackage{soul,color}
\newtheorem{theorem}{Theorem}[section]
\newtheorem{proposition}[theorem]{Proposition}

\theoremstyle{remark}
\newtheorem{remark}[theorem]{Remark}
\theoremstyle{example}

\numberwithin{equation}{section}

\begin{document}

\newcommand{\Surf}{{\partial \Omega}}
\newcommand{\Vol}{{\Omega}}
\newcommand{\dVol}{{\Gamma}}
\newcommand\txtred[1]{{\color{red}#1}}
\newcommand\txtblue[1]{{\color{blue}#1}}
\newcommand{\xj}{{x_j}}
\newcommand{\ui}{{u_i}}
\newcommand{\uj}{{u_j}}
\newcommand{\pd}[2]{\frac{\partial #1}{\partial #2}}
\newcommand{\taupr}{\tau^\prime}
\newcommand{\Rot}{\mathbb{R}}

\newcommand{\half}{\frac{1}{2}}
\newcommand{\BigO}[1]{\ensuremath{\operatorname{O}\bigl(#1\bigr)}}
\newcommand{\Order}{\mathcal{O}}
\newcommand{\U}{{\mathbf{U}}}
\newcommand{\Ue}{{\mathbf{V}}}
\newcommand{\Ux}{{\mathbf{U}_x}}
\newcommand{\Vv}{{\mathbf{V}}}
\newcommand{\vv}{{\overline{v}}}
\newcommand{\G}{{\mathbf{G}}}
\newcommand{\At}{\tilde{A}}
\newcommand{\Pp}{\mathcal{P}}
\newcommand{\Ppo}{{\mathcal{\overline{P}}}}
\newcommand{\Q}{\mathcal{Q}}
\newcommand{\Qx}{\mathcal{Q}_{x}}
\newcommand{\Qxe}{\mathcal{Q}_x^e}
\newcommand{\Qxx}{\mathcal{Q}_{xx}}
\newcommand{\Qxi}{\mathcal{Q}_{x^i}}
\newcommand{\B}{\mathcal{B}}
\newcommand{\Lagr}{\mathcal{L}}
\newcommand{\lagr}{\mathcal{l}}
\newcommand{\D}{\mathcal{D}}
\newcommand{\Dx}{\mathcal{D}_{x}}
\newcommand{\Dxi}{\mathcal{D}_{x^i}}
\newcommand{\Dxa}{\mathcal{D}_{x^1}}
\newcommand{\Dxb}{\mathcal{D}_{x^2}}
\newcommand{\Dxc}{\mathcal{D}_{x^3}}
\newcommand{\Bx}{\B\Dx}
\newcommand{\Lx}{\mathbb{L}_x}
\newcommand{\Di}{\mathcal{D}_{AD}}
\newcommand{\Dii}{\mathcal{D}_{AD_i}}
\newcommand{\De}{\mathcal{D}^{e}}
\newcommand{\Dxx}{\mathcal{D}_{xx}}
\newcommand{\Ee}{\mathcal{E}}
\newcommand{\R}{\mathcal{R}}
\newcommand{\nn}{\nonumber}
\newcommand{\Nn}{\ell}
\newcommand{\Nns}{\ell^\star}
\newcommand{\Mm}{\mathcal{M}}
\newcommand\norm[1]{\left\lVert#1\right\rVert}
\newcommand\bint[3]{\left. #1 \right|_{#2}^{#3}}
\newcommand\sg[1]{\sigma^#1}
\newcommand\eps{\varepsilon}
\newcommand\epst{\overline{\varepsilon}}
\newcommand{\V}{\overline{u}}
\newcommand{\F}{\mathcal{F}}
\newcommand{\Aa}{\mathcal{A}}
\newcommand\bb[1]{\textcolor{blue}{#1}}
\newcommand\rred[1]{\textcolor{red}{#1}}
\newcommand\ggreen[1]{\textcolor{green}{#1}}
\newcommand\bblue[1]{\textcolor{blue}{#1}}
\newcommand{\I}{\mathcal{\tilde{I}}}
\newcommand{\A}{\mathcal{\tilde{A}}}
\newcommand{\Iha}{\mathcal{\hat{I}}_{4t2}}
\newcommand{\Ihb}{\mathcal{\hat{I}}_{4t4}}
\newcommand{\al}{\alpha}
\newcommand\bft[1]{{\bf{#1}}}
\newcommand{\drho}{\Delta\rho}
\newcommand{\BT}{\mathbb{BT}}
\newcommand{\BTI}{\mathbb{BTI}}
\newcommand{\SAT}{\mathbb{SAT}}
\newcommand{\BC}{\mathbb{BC}}
\newcommand{\Ut}{\tilde{U}}
\newcommand{\C}{\mathbb{C}}

\makeatletter
\newcommand{\doublewidetilde}[1]{{%
  \mathpalette\double@widetilde{#1}%
}}
\newcommand{\double@widetilde}[2]{%
  \sbox\z@{$\m@th#1\widetilde{#2}$}%
  \ht\z@=.9\ht\z@
  \widetilde{\box\z@}%
}
\makeatother


\begin{frontmatter}

%

\title{A well posed and stable canonical evaporation model problem for phase-change in two-phase flows}

\author[sweden,southafrica]{Jan Nordstr\"{o}m}
\corref{firstcorrespondingauthor}
\cortext[firstcorrespondingauthor]{Corresponding author}
\ead{jan.nordstrom@liu.se}
\address[sweden]{Department of Mathematics,  Applied Mathematics, Link\"{o}ping University, SE-581 83 Link\"{o}ping, Sweden}
\address[southafrica]{Department of Mathematics   $\&$ Applied Mathematics, University of Johannesburg, Auckland Park 2006, South Africa}

\begin{abstract}
We formulate a well posed interface formulation for canonical one-dimensional evaporation two-phase model problems (the Stefan and Sucking problems) commonly used to  validate  
production codes. 
We focus on the interface between the vapor and the liquid and derive conditions leading to an energy bound and well-posedness. Next, by mimicking the continuous analysis, we discretize using high order accurate numerical methods on summation-by-parts form, impose the interface conditions weakly and prove energy stability.
\end{abstract}



\end{frontmatter}


\section{Introduction}



We consider two 
model problems for incompressible two-phase liquid-gas flows in the volume-of-fluid formulation  \cite{Hirt1981}.
In particular, we focus on the evaporation interface for 
the so called Stefan and Sucking problems \cite{ph-c}.
These extremely stiff model problems are commonly used to  validate  
production codes \cite{Flow3D}. 
We provide a detailed analysis of the given interface model used in these problems.



After transformation of the moving evaporation interface problems to a steady frame 
\cite{NIKKAR201582,move1}
we write the equations in skew-symmetric form 
\cite{NORDSTROM2024_BC,nordstrom2024skewsymmetric_jcp} 
suitable 
for numerical differentiation 
and derive energy bounding well posed interface conditions. 
Following \cite{nordstrom_roadmap}, we discretize using 
numerical methods on Summation-By-Parts (SBP) form, impose the interface conditions weakly using the Simultaneous Approximation Technique (SAT)  \cite{svard2014review,fernandez2014review} 
and prove energy stability. 

 \section{The initial boundary value problem (IBVP)}
\label{sec:CF}
The 
governing coupled IBVP
in the domain $\Omega = [x_0,x_n]$  describing the interface action at $x_\delta(t)$ 
are
\begin{subequations}
\label{eq:problem_orig}
\begin{equation}
\label{eq:problemv}
\begin{split}
    \beta_v\Big(\big(T_v\big)_t + u_v \big(T_v\big)_x\Big) &= k_v \big(T_v\big)_{xx}, \;\;\; x \in [x_0,x_\delta],\\
    T_v &= T_\delta, \;\;\; x = x_\delta(t)
\end{split}    
\end{equation}
\begin{equation}
\label{eq:probleml}
\begin{split}
    \beta_l\Big(\big(T_l\big)_t + u_l \big(T_l\big)_x \Big) &= k_l \big(T_l\big)_{xx}, \;\;\; x \in [x_\delta, x_n],\\
     T_l &= T_\delta, \;\;\; x = x_\delta(t)
\end{split}    
\end{equation}
\end{subequations}
where we ignore the outer boundary and initial conditions.
In (\ref{eq:problemv}) and (\ref{eq:probleml}): $T, \rho,  C_p, k$  denotes temperature, density, specific heat at constant pressure and scaled heat conduction coefficient respectively while $\beta= \big(\rho C_p\big)$. The subscripts $v,l,\delta$ denotes vapour, liquid and evaporation interface respectively. 
The problems are connected at the time-dependent moving interface $x_\delta(t)$ by the externally given evaporation temperature $T_\delta$. 
We also assume that the wave speeds $u_v$ and $u_l$ are positive and satisfy the incompressibility constraint.  
\begin{remark}
Stable outer boundary conditions for the advection-diffusion problem (\ref{eq:problem_orig}) are given  
in \cite{Carpenter1999341}.
\end{remark}

To transform the coupled problems to a steady frame, we introduce the coordinate transformations \cite{NIKKAR201582,move1}
$x = x(\xi, \tau ), t = \tau, \xi = \xi(x, t)$ for $ x \in [x_0,x_\delta]$ and 
$x = x(\eta, \tau ), t = \tau, \eta = \eta(x, t)$ for $x \in [x_\delta,x_n]$.
These transformations 
lead to the metric identities and geometric conservation laws \cite{lesoinne1996geometric}
\begin{equation}
\label{eq:transformation-metric}
\begin{split}
J_v \xi_x = 1  \;\;\; & J_v \xi_t = -x_\tau, \;\;\;  \text{and} \;\;\; (J_v)_\tau + (J_v \xi_t)_\xi=0, \;\;\;\; x \in [x_0,x_\delta],  \\
J_l \eta_x = 1  \;\;\; & J_l \eta_t = -x_\tau, \;\;\;  \text{and} \;\;\; (J_l)_\tau + (J_l \eta_t)_\eta=0, \;\;\;\; x \in [x_\delta, x_n].
\end{split}
\end{equation}
The relations  (\ref{eq:transformation-metric}) 
transform equations (\ref{eq:problemv}) and (\ref{eq:probleml}) on moving domains to equations on fixed domains
\begin{subequations}
    \label{eq:main}
\begin{equation}\label{eq:mainv}
\begin{split}
    \beta_v\Big(\big(J_v T_v\big)_\tau  + \big(a_vT_v\big)_\xi\Big) &= k_v\big(J^{-1}_v T_{v,\xi}\big)_\xi, \;\;\; \xi\in\Omega_v = [\xi_0,\delta],\\
    T_v& = T_\delta, \;\;\; \xi = \delta
\end{split}
\end{equation}
\begin{equation}\label{eq:mainl}
\begin{split}
    \beta_l\Big(\big(J_l T_l\big)_\tau  + \big(a_lT_l\big)_\eta\Big) &= k_l\big(J^{-1}_l T_{l,\eta}\big)_\eta, \;\;\; \eta\in\Omega_l = [\delta,\eta_n],\\
    T_l &= T_\delta, \;\;\; \eta = \delta
\end{split}
\end{equation}
\end{subequations}
where  $\partial T_v/\partial \xi =T_{v,\xi}$, $\partial T_l/\partial \eta =T_{l,\eta}$.  $J_v$ and $J_l$ are the Jacobians  and the transformed wave speeds are
\begin{equation}
    \label{eq:a-devs}
a_v=u_v+J_v \xi_t=u_v -x_\tau=u_v-\Tilde{u}_v \quad \text{and} \quad  a_l=u_l+J_v \eta_t=u_l -x_\tau=u_l-\Tilde{u}_l.
\end{equation}

At the interface (indicated with subscript $\delta$) the mass flow rate is given by the relation
\begin{equation}
\label{eq:conservation_int}
    \rho_v\big((u_v)_\delta - \Tilde{u}_\delta\big) = \rho_l \big((u_l)_\delta - \Tilde{u}_\delta\big), \;\;\;\text{or} \;\;\; (u_l)_\delta = \gamma (u_v)_\delta + \big( 1 - \gamma \big) \Tilde{u}_\delta
\end{equation}
where $\gamma = \rho_v / \rho_l$ and the mesh velocity $\Tilde{u}_\delta$ at the interface $x_\delta$ 
is given by the jump in temperature flux as 
\begin{equation}
    \label{eq:meshvelocity}
    \rho_v h_{lv} \Tilde{u}_\delta = \big[k_l J_l^{-1} T_{l,\eta} - k_v J_v^{-1} T_{v, \xi} \big]_\delta.
\end{equation}
Furthermore, \eqref{eq:conservation_int}  
gives us the relations between the transformed wave speeds at the interface $x=x_\delta$  as
\begin{equation}
\label{eq:af}
        (a_v)_\delta = (- x_\tau  + u_v)_\delta = (u_v)_\delta-\Tilde{u}_\delta \quad and \quad  (a_l)_\delta = (-x_\tau + u_l)_\delta = \gamma ((u_v)_\delta-\Tilde{u}_\delta)=\gamma  (a_v)_\delta.
\end{equation}
\begin{remark}
The time-integral of interface velocity 
is bounded 
since $\int_0^{t} \Tilde{u}_\delta  \ d\tau  = \int_0^{t}  x_\tau   \ d\tau  = x(t)-x(0) < \Omega$.
\end{remark}
\begin{remark}
The linear IBVPs   (\ref{eq:problem_orig}) are now transformed to a quasi-nonlinear setting in (\ref{eq:main}) since $a_v$ and $a_l$ are solution dependent
via the varying metric coefficients which in turn vary with the mesh velocity (\ref{eq:meshvelocity}).
\end{remark}

\section{Continuous energy analysis}
Applying the energy method (multiply with the solution and integrate over the domain) to \eqref{eq:main} 
yields
\begin{equation}
\label{eq:semi_energy}
    \dfrac{d}{d\tau}\Big(\beta_v\|T_v\|^2_{J_v} + \beta_l\|T_l\|^2_{J_l} \Big)  + 2 \Big (k_v  \|J_v^{-1}T_{v,\xi}\|^2_{J_v} +k_l  \|J_l^{-1}T_{l,\eta}\|^2_{J_l}\Big) = IT + BT,
\end{equation}
where $\|T_v\|^2_{J_v} =\int_0^\delta T_v^2 J_v d\xi$ and $\|T_l\|^2_{J_l} =\int_\delta^{\eta_n}\ T_l^2 J_l d\eta $.  The interface term denoted $IT$ is given by 
\begin{equation}
\label{eq:IBT}
        IT = \big(-a_v\beta_v T_v^2 + 2 k_v T_v J_v^{-1} T_{v,\xi}\big)_\delta  +\big(a_l\beta_lT_l^2 - 2 k_l T_l J_l^{-1} T_{l,\eta}\big)_\delta.
\end{equation}
As stated above, we ignore the fixed outer boundary term denoted $BT$ and focus on evaporation interface.

The interface condition $T_v=T_l=T_\delta$ and mesh velocity (\ref{eq:meshvelocity}) together with the relations \eqref{eq:af} yields
\begin{equation}
\label{eq:IT1}
    \begin{split}
         IT & = 
         (\beta_l a_l-\beta_v a_v)_\delta T_\delta^2 - 2 T_\delta \big(k_l J_l^{-1} T_{l,\eta} - k_v J_v^{-1} T_{v,\xi})_\delta\\
         &= \big(\beta_l \gamma - \beta_v\big)(a_v)_\delta T_\delta^2  - 2  \big(T_\delta \rho_v h_{lv}\big) \Tilde{u}_\delta \\
         &= C_0 \big((u_v)_\delta - \Tilde{u}_\delta\big)  - C_1 \Tilde{u}_\delta,
    \end{split}
\end{equation}
where $C_0 = (\beta_l \gamma - \beta_v) T_\delta^2= \rho_v\big[\big(C_P\big)_l - \big(C_p\big)_v\big]T_\delta^2>0 $ and $C_1 = 2 \big(T_\delta \rho_v h_{lv}\big) >0$ are positive constants. 

We have proved Proposition (\ref{lemma:strong_cont}).
\begin{proposition}\label{lemma:strong_cont}
The strong interface treatment for equations (\ref{eq:main}) with the mesh velocity definitions (\ref{eq:a-devs})-\eqref{eq:af} 
is dissipative and energy bounded if  $(a_v)_\delta <  0$ and  $\Tilde{u}_\delta \ge 0$.  It is energy bounded otherwise.
\end{proposition}

\subsection{A continuous weak imposition of the interface conditions}
The relations (\ref{eq:main}) with a weak impositions of the interface condition $T_v=T_l=T_\delta$  becomes
\begin{equation}\label{eq:main+penn}
\begin{split}
    \beta_v\Big(\big(J_v T_v\big)_\tau  + \big(a_vT_v\big)_\xi\Big) &= k_v\big(J^{-1}_v T_{v,\xi}\big)_\xi + L_\delta\Big(\sigma_v^1(T_v-T_\delta) + \sigma_v^2(T_v-T_l) + \sigma_v^3 (\partial /\partial \xi)^T(T_v-T_\delta)\Big)\\
    \beta_l\Big(\big(J_l T_v\big)_\tau  + \big(a_lT_l\big)_\eta\Big) &= k_l\big(J^{-1}_l T_{v,\xi}\big)_\eta + L_\delta\Big(\sigma_l^1(T_l-T_\delta) + \sigma_l^2(T_l-T_v) + \sigma_l^3 (\partial /\partial \eta)^T(T_l-T_\delta)\Big).
\end{split}
\end{equation}
In (\ref{eq:main+penn}),  $L_\delta$ is a lifting operator \cite{Arnold20011749} imposing the interface conditions at $x=x_\delta$ weakly. The lifting operator is defined by $ \int  \psi L_\delta (\phi) dx = (\psi \phi)_\delta $ for two smooth functions $\psi$ and $\phi$. The six penalty coefficients $\sigma_k^i$ will be determined below for stability. We also use the non-conventional transposed derivative operator $(\partial /\partial x)^T$ defined such that for a smooth function
$\phi$, multiplication from the left yields $\phi(\partial /\partial x)^T=\phi_x$.
\begin{remark}
Existence and well-posedness require a minimum number of boundary and interface conditions \cite{nordstrom2020}.
The coupled problem (\ref{eq:main+penn}) 
include only two interface conditions
since $T_v \equiv T_l$ is trivially true.
\end{remark}

By again applying the energy method as above, and ignoring the outer boundary terms $BT$  we obtain 
\begin{equation}
\label{eq:semi_energy_SAT}
    \dfrac{d}{d\tau}\Big(\beta_v\|T_v\|^2_{J_v} + \beta_l\|T_l\|^2_{J_l} \Big)  + 2 k_v  \|J_v^{-1}T_{v,\xi}\|^2_{J_v} +2 k_l  \|J_l^{-1}T_{l,\eta}\|^2_{J_l} = IT + SAT
\end{equation}
where $IT$ is given in (\ref{eq:IBT}) and $SAT=SAT_v+SAT_l$ is the additional weak contribution given by 
\begin{equation}\label{eq:SAT}
\begin{split}
SAT_v &= 2\Big(\sigma_v^1T_v(T_v-T_\delta) + \sigma_v^2T_v(T_v-T_l) + \sigma_v^3 T_{v,\xi} (T_v-T_\delta)\Big)_\delta,\\
SAT_l &= 2\Big(\sigma_l^1T_l(T_l-T_\delta) + \sigma_l^2T_l(T_l-T_v) + \sigma_l^3 T_{l,\eta} (T_l-T_\delta)\Big)_\delta.
\end{split}
\end{equation}

First we scale the gradient terms in $SAT$ by letting $\sigma_v^3 = \sigma_v^4 k_v  J_v^{-1}$ and $\sigma_l^3 = \sigma_l^4 k_l  J_l^{-1}$ to make them comparable to the $IT$ ones. Next, we group similar terms from $IT$ and $SAT$  denoting them GRAD to get
\begin{equation}\label{eq:GRAD1}
GRAD= 2\Big((\sigma_v^4+1) k_v  J_v^{-1} T_v T_{v,\xi} + (\sigma_l^4-1) k_l  J_v^{-1} T_l T_{l,\eta} - \sigma_v^4k_v  J_v^{-1} T_\delta T_{v,\xi}- \sigma_l^4k_l  J_l^{-1} T_\delta T_{l,\eta}\Big)_\delta.
\end{equation}
An obvious first choice is to let $\sigma_v^4 = -1$ and $\sigma_l^4 = +1$ which remove the first two term and leaves us with
\begin{equation}\label{eq:GRAD2}
GRAD= 2 T_\delta \big(k_v  J_v^{-1} T_{v,\xi}- k_l  J_l^{-1} T_{l,\eta}\big)_\delta = - 2 \big(T_\delta \rho_v h_{lv}\big)\Tilde{u}_\delta 
\end{equation}
by using definition 
(\ref{eq:meshvelocity}). The coefficient multiplying $\Tilde{u}_\delta$ is given by external data, it is negative and constant. Hence, the gradient terms provide dissipation if $\Tilde{u}_\delta >0$ and bounded energy growth if $\Tilde{u}_\delta <0$.

Next we consider the terms consisting of point values.  Grouped together we denote them $PVAL$ where
\begin{equation}\label{eq:PVAL1}
PVAL= \begin{pmatrix}  T_v \\ T_l  \\  T_\delta \end{pmatrix}^T
\begin{pmatrix} -\beta_v a_v + 2 (\sigma_v^1 +\sigma_v^2)& -(\sigma_v^2 +\sigma_l^2)                          & -\sigma_v^1 \\  
                          -(\sigma_v^2 +\sigma_l^2)                          &  +\beta_l a_l + 2 (\sigma_l^1 +\sigma_l^2) & -\sigma_l^1 \\ 
                          -\sigma_v^1                                                 & -\sigma_l^1                                                 & 0
\end{pmatrix}                           
 \begin{pmatrix}  T_v \\ T_l  \\  T_\delta \end{pmatrix}.
\end{equation}
Setting $\sigma_v^2 = \sigma_l^2=-\sigma/2$ give  $PVAL=PVAL_2+ PVAL_v+PVAL_l$ where  $PVAL_2=-\sigma(T_v-T_l)^2$ and
\begin{equation}\label{eq:PVAL2}
PVAL_v= \begin{pmatrix}  T_v \\ T_\delta \end{pmatrix}^T
\begin{pmatrix} -\beta_v a_v + 2 \sigma_v^1 & -\sigma_v^1 \\  
                                                   -\sigma_v^1 & 0
\end{pmatrix}                         
 \begin{pmatrix}  T_v \\ T_\delta \end{pmatrix}, \ 
 PVAL_l= \begin{pmatrix}  T_l \\ T_\delta \end{pmatrix}^T
\begin{pmatrix}  +\beta_l a_l + 2 \sigma_l^1 & -\sigma_l^1 \\  
                                                   -\sigma_l^1 & 0
\end{pmatrix}                         
 \begin{pmatrix}  T_v \\ T_\delta \end{pmatrix}.
\end{equation}
Now there are two options depending on the sign of $(a_v)_\delta$ (recall that $(a_l)_\delta=\gamma (a_v)_\delta$). 

Firstly, for  $(a_v)_\delta>0$,  
then let $\sigma_v^1=0$ and $\sigma_l^1=-\beta_l (a_l)_\delta$. This leads to
\begin{equation}\label{eq:PVAL3}
PVAL=(\beta_l a_l-\beta_v a_v)_\delta T_\delta^2-\sigma(T_v-T_l)^2-\beta_l (a_l)_\delta(T_l-T_\delta)^2 -\beta_v (a_v)_\delta (T_v^2-T_\delta^2).
\end{equation}

Secondly, for $(a_v)_\delta<0$, 
then let $\sigma_l^1=0$ and $\sigma_v^1=\beta_v (a_v)_\delta $. This gives
\begin{equation}\label{eq:PVAL5}
PVAL=(\beta_l a_l-\beta_v a_v)_\delta T_\delta^2 -\sigma(T_v-T_l)^2+\beta_v (a_v)_\delta(T_v-T_\delta)^2 +\beta_l (a_l)_\delta (T_l^2-T_\delta^2).
\end{equation}
The $1st$ terms in (\ref{eq:PVAL3}),(\ref{eq:PVAL5}) is the strong result, the $2nd$ and $3rd$ one provide dissipation and the $4th$ one is bounded.
Combining GRAD in (\ref{eq:GRAD2}) with PVAL in (\ref{eq:PVAL3}) for $(a_v)_\delta <  0$ and with (\ref{eq:PVAL5}) $(a_v)_\delta >  0$ yields
\begin{equation}\label{eq:IT+SAT-weak-strong-avneg}
IT + SAT=(\beta_l a_l-\beta_v a_v)_\delta T_\delta^2 - 2 \big(T_\delta \rho_v h_{lv}\big)\Tilde{u}_\delta -\sigma(T_v-T_l)^2+\beta_v (a_v)_\delta(T_v-T_\delta)^2 +\beta_l (a_l)_\delta (T_l^2-T_\delta^2),
\end{equation}
\begin{equation}\label{eq:IT+SAT-weak-strong-avpos}
IT + SAT=(\beta_l a_l-\beta_v a_v)_\delta T_\delta^2-2 \big(T_\delta \rho_v h_{lv}\big)\Tilde{u}_\delta -\sigma(T_v-T_l)^2-\beta_l (a_l)_\delta(T_l-T_\delta)^2 -\beta_v (a_v)_\delta (T_v^2-T_\delta^2)
\end{equation}
respectively.
Both (\ref{eq:IT+SAT-weak-strong-avneg}) and (\ref{eq:IT+SAT-weak-strong-avpos}) leads to energy stability via time-integration of the energy rate (\ref{eq:semi_energy_SAT}).

We have now proved Proposition (\ref{lemma:weak_cont}).
\begin{proposition}\label{lemma:weak_cont}
Consider the weak interface treatment formulation (\ref{eq:main+penn}) for equations (\ref{eq:main}) with the mesh velocity definitions (\ref{eq:a-devs})-\eqref{eq:af}.  Furthermore, let the penalty coefficients 
 be
\begin{equation}\label{eq:pen-weak-strong-avneg}
(a_v)_\delta <  0:  \sigma_v^1=\beta_v (a_v)_\delta, \quad \sigma_l^1=0,  \quad \sigma_v^2=\sigma_l^2=-\sigma/2 \leq 0, \quad \sigma_v^3 = -k_v  J_v^{-1}, \quad \sigma_l^3 = k_l  J_l^{-1}
\end{equation}
\begin{equation}\label{eq:pen-weak-strong-avpos}
(a_v)_\delta <  0: \sigma_v^1=0, \quad \sigma_l^1=-\beta_l (a_l)_\delta,  \quad \sigma_v^2=\sigma_l^2=-\sigma/2 \leq 0, \quad \sigma_v^3 = -k_v  J_v^{-1}, \quad \sigma_l^3 = k_l  J_l^{-1}.
\end{equation}
Then the formulation (\ref{eq:main+penn}) is energy bounded. For $(a_v)_\delta=0$ we simply let $\sigma_v^1=\sigma_l^1=0$. 
\end{proposition}
 
\section{Semi-discrete energy analysis}
Following \cite{nordstrom_roadmap},
we mimic the continuous analysis numerically by using the SBP-SAT framework. 
We denote $\mathbf{T}_v$ as the numerical approximation of $T_v$ at all grid points 
and the vectors $D_\xi \mathbf{T}_v=\mathbf{T}_{v,\xi} $  to be the approximation of $T_{v,\xi}$. $D_\xi = P^{-1}_\xi Q_\xi$ is the differention matrix on SBP form \cite{svard2014review} with the SBP property
\begin{equation}\label{eq:SBP}
	Q_\xi + Q^T_\xi = B_\xi =  E_N^\xi - E_0^\xi, \quad where \quad  E_0^\xi  = \rm{diag}(1,0,\dots,0) \quad and  \quad E_N^\xi  = \rm{diag}(0,0,\dots,1).
\end{equation}
The matrix $P_\xi$ is diagonal, positive definite and defines a quadrature formula ${\boldsymbol \phi}^T P_\xi {\boldsymbol \psi} \approx  \int \phi \psi dx $ and a norm  ${\boldsymbol \phi}^T P_\xi {\boldsymbol \phi} =  ||{\boldsymbol \phi} ||^2_{P_\xi}$ where $ {\boldsymbol \phi}$,${\boldsymbol \psi}$ are vectors with  $\phi_i$,$\psi_i$  at node $\xi_i$. The same definitions hold in the $\eta$ direction.

\subsection{A semi-discrete weak imposition of the interface conditions}
The semi-discrete version of the continuous equations (\ref{eq:main+penn}) with SAT interface conditions are
 \begin{equation}
        \label{eq:mainvdis_tot}
        \begin{split}
            \dfrac{\beta_v }{2}\Big(\big(\mathbf{J}_v \mathbf{T}_v\big)_\tau + \mathbf{J}_v \mathbf{T}_{v,\tau}  + D_\xi A_v \mathbf{T}_v+ A_v D_\xi \mathbf{T}_v \Big)  &= k_v D_{\xi} \big(\mathbf{J}_v^{-1}\mathbf{T}_{v,\xi} \big)\\   
             + P_\xi^{-1}\Big(\sigma_v^1E^\xi_N (\mathbf{T}_v-\mathbf{T}^\delta_v) &+ \sigma_v^2 E^\xi_N (\mathbf{T}_v-\mathbf{T}^l_v) + \sigma_v^3 D_\xi ^TE^\xi_N (\mathbf{T}_v-\mathbf{T}^\delta_v)\Big),\\                                                                                                                                                                                                                                                                                   
            \dfrac{\beta_l }{2}\Big(\big(\mathbf{J}_l \mathbf{T}_l\big)_\tau + \mathbf{J}_l \mathbf{T}_{l,\tau}  + D_\eta A_l \mathbf{T}_l+ A_l D_\eta \mathbf{T}_l \Big) & = k_l D_{\eta} \big(\mathbf{J}_l^{-1}\mathbf{T}_{l,\eta} \big)\\   
             + P_\eta^{-1}\Big(\sigma_l^1 E^\eta_0 (\mathbf{T}_l-\mathbf{T}^\delta_l) &+ \sigma_l^2E^\eta_0 (\mathbf{T}_l-\mathbf{T}^v_l) + \sigma_l^3 D_\eta^T E^\eta_0 (\mathbf{T}_l-\mathbf{T}^\delta_l) \Big).                                                                                                                                                                                                                                                                                        
        \end{split}
    \end{equation}
In (\ref{eq:mainvdis_tot}), $\mathbf{T}^\delta_v$ and $\mathbf{T}^l_v$  are equal to $\mathbf{T}_v$ except in the last position which contain $T_\delta$ and $(\mathbf{T}_l)_0$. Similary   $\mathbf{T}^\delta_l$ and $\mathbf{T}^v_l$ are equal to $\mathbf{T}_l$ except in the first position which contain $T_\delta$ and $(\mathbf{T}_v)_N$. Furthermore, $\mathbf{J}_v, \mathbf{J}_l$ and $A_v, A_l$ are diagonal matrices 
approximating $J_v,J_l$ and $a_v, a_l$  while $P_\xi^{-1}$, $P_\eta^{-1}$  are the discrete lifting operators.

The energy method applied to the two equations (\ref{eq:mainvdis_tot}) and adding the results leads to
 \begin{equation}
        \label{eq:mainvdis_tot_energy}
        \dfrac{d}{d\tau} \big(\beta_v \|\mathbf{T}_v\|^2_{J_vP_\xi}  + \beta_l \|\mathbf{T}_l\|^2_{J_lP_\eta}\big) + 2 \big( k_v \|\mathbf{J}_v^{-1}\mathbf{T}_{v,\xi}\|^2_{J_vP_\xi} + k_l \|\mathbf{J}_l^{-1}\mathbf{T}_{l,\eta}\|^2_{J_lP_\eta} \big) = IT+SAT
    \end{equation}
where we used  
(\ref{eq:SBP}).
The semi-discrete interface term $IT$ and 
$SAT=SAT_v+SAT_l$  term are given below
\begin{equation}
\label{eq:IBT_d}
        IT = -\beta_v\mathbf{T}_v^T E^\xi_N A_v \mathbf{T}_v + 2k_v \mathbf{T}_v^T E^\xi_N \mathbf{J}_v^{-1} \mathbf{T}_{v,\xi} +\beta_l\mathbf{T}_l^T E^\eta_0 A_l \mathbf{T}_l - 2k_l \mathbf{T}_l^T E^\eta_0  \mathbf{J}_l^{-1} \mathbf{T}_{l,\eta},
\end{equation}
 \begin{equation}
        \label{eq:SAT_tot}
        \begin{split}
           SAT_v &=  
                           2\Big(\sigma_v^1(E^\xi_N \mathbf{T}_v)^T (\mathbf{T}_v-\mathbf{T}^\delta_v) + \sigma_v^2(E^\xi_N \mathbf{T}_v)^T (\mathbf{T}_v-\mathbf{T}^l_v) + \sigma_v^3 (E^\xi_N \mathbf{T}_{v,\xi})^T (\mathbf{T}_v-\mathbf{T}^\delta_v)\Big),\\    
          SAT_l &=  
                           2\Big(\sigma_l^1(E^\eta_0 \mathbf{T}_l)^T (\mathbf{T}_l-\mathbf{T}^\delta_v) + \sigma_l^2(E^\eta_0 \mathbf{T}_l)^T (\mathbf{T}_l-\mathbf{T}^v_l) + \sigma_l^3 (E^\eta_0 \mathbf{T}_{l,\eta})^T (\mathbf{T}_l-\mathbf{T}^\delta_l)\Big).                                                                                                                                                                                                                                                                                                                                                                                                                                                                                                                                                              
        \end{split}
    \end{equation}
The semi-discrete $IT$  (\ref{eq:IBT_d}) correspond term by term to the continuous $IT$ (\ref{eq:IBT}) where 
$E^\xi_N$ and $E^\eta_0$  pick out only the relevant interface terms. 
Hence (\ref{eq:IBT}) holds also in the semi-discrete case.  The semi-discrete $SAT=SAT_v+SAT_l$  also correspond term by term to the continuous  $SAT$ in (\ref{eq:SAT}).  The difference is that gradients are computed numerically.
This implies that the continuous 
derivation and 
Proposition (\ref{lemma:weak_cont}) holds also in the semi-discrete case (with the differences mentioned above) which proves Proposition (\ref{lemma:weak_disc}).
\begin{proposition}\label{lemma:weak_disc}
Consider the semi-discrete 
formulation (\ref{eq:mainvdis_tot}) for equations (\ref{eq:main}) with the mesh velocity definitions (\ref{eq:a-devs})-\eqref{eq:af} and gradients computed numerically.
Furthermore, let the penalty coefficients 
be
\begin{equation}\label{eq:pen-weak-strong-avneg_d}
(a_v)_\delta <  0:  \sigma_v^1=\beta_v (a_v)_\delta, \quad \sigma_l^1=0,  \quad \sigma_v^2=\sigma_l^2=-\sigma/2 \leq 0, \quad \sigma_v^3 = -k_v  J_v^{-1}, \quad \sigma_l^3 = k_l  J_l^{-1}
\end{equation}
\begin{equation}\label{eq:pen-weak-strong-avpos_d}
(a_v)_\delta <  0: \sigma_v^1=0, \quad \sigma_l^1=-\beta_l (a_l)_\delta,  \quad \sigma_v^2=\sigma_l^2=-\sigma/2 \leq 0, \quad \sigma_v^3 = -k_v  J_v^{-1}, \quad \sigma_l^3 = k_l  J_l^{-1}.
\end{equation}
Then the formulation (\ref{eq:mainvdis_tot}) is energy stable. For $(a_v)_\delta=0$ we simply let $\sigma_v^1=\sigma_l^1=0$. 
\end{proposition}

\section{Summary}\label{sec:conclusion}
We formulated a provably energy bounded interface formulation for two canonical one-dimensional evaporation two-phase model problems, the so called Stefan and Sucking problems. These extremely stiff model problems are commonly used to validate existing production codes.
Based on the energy bounded continuous formulation, we discretized using numerical methods 
on SBP-SAT form 
and proved energy stability.  

\section*{Acknowledgments}
\vspace{-0.1truecm}
\noindent J. N. was supported by 
University of Johannesburg Global Excellence and Stature Initiative Funding. 
J.N. acknowledge discussions with A. G. Malan and H. Hanif from the InCFD Group at University of Cape Town.

\bibliographystyle{elsarticle-num}
\bibliography{References_Jan}

\end{document}